# GAUSSIAN PERTURBATIONS OF CIRCLE MAPS: A SPECTRAL APPROACH

### By John Mayberry

*Cornell University*

In this work, we examine spectral properties of Markov transition operators corresponding to Gaussian perturbations of discrete time dynamical systems on the circle. We develop a method for calculating asymptotic expressions for eigenvalues (in the zero noise limit) and show that changes to the number or period of stable orbits for the deterministic system correspond to changes in the number of limiting modulus 1 eigenvalues of the transition operator for the perturbed process. We call this phenomenon a $\lambda$-bifurcation. Asymptotic expressions for the corresponding eigenfunctions and eigenmeasures are also derived and are related to Hermite functions.

**1. Introduction.** Studies of mechanical and biological oscillators have suggested that the eigenvalues of Markov transition operators can be used to analyze bifurcation behavior in random perturbations of deterministic systems (see [3, 9, 10, 11]). In particular, numerical observations in these papers show that bifurcations in the underlying deterministic system often correspond to changes in the number of eigenvalues with modulus close to 1 when the perturbation size is small. In this paper, we follow up on these numerical observations by providing a rigorous example in which this phenomenon occurs.

To this end, we consider Gaussian perturbations of dynamical systems on the circle exhibiting stable periodic behavior. We will provide a method for calculating limiting eigenvalues (as $\varepsilon \to 0$) of the transition operator for the system

$$X_{n+1}^{\varepsilon} = f(X_n^{\varepsilon}) + \varepsilon \sigma(X_n^{\varepsilon})\chi_n \bmod 2\pi, \tag{1}$$

where $f$ is a sufficiently smooth circle map with a finite number of periodic orbits that attract all other orbits of $x_{n+1} = f(x_n)$, $\{\chi_n\}_{n=0}^{\infty}$ is a family of i.i.d. standard normal random variables and $\sigma$ is sufficiently smooth and positive. Our main results are described in Section 2 and basically state that the limiting eigenvalues of $T^\varepsilon$ are determined by the derivative of $f$ along periodic orbits while the corresponding limiting eigenvectors are related to hermite functions. To illustrate our methods, Sections 3–5 give a detailed analysis of the case when $f$ has one stable fixed point $x_s$ and one unstable fixed point $x_u$. Section 3 describes the setup and basic results in this setting and develops our primary tool: a block decomposition of the transition operator which allows us to calculate spectral properties by focusing on the "local" action of the transition operator near the fixed points of $f$. Sections 4 and 5 then contain details of the local analysis and extensions to general periodic orbits are derived in Section 6.

Before moving on to the body of our work, we note that our basic model (1) can be used as an approximate heuristic for studying the dynamics of sequences of firing phases in integrate-and-fire models with a white noise component which provides some connection between our results and the numerical observations in [9, 10, 11]. This connection will be further developed in the paper [1]. We also leave it to the reader to check that many of our results concerning eigenvalues of $T^\varepsilon$ remain true if we replace the $\chi_n$ in (1) with some other sequence of i.i.d., finite moment generating function random variables. We focus here on the Gaussian case since the calculations then yield particularly interesting formulas for eigenvectors in terms of hermite functions. An interesting question for future research would be extensions to the case when the asymptotic behavior of the deterministic system $x_{n+1} = f(x_n)$ is chaotic although this is likely to require different techniques (see, e.g., [7]).

## 2. General heuristic.

In this section, we describe our setting and main results. Throughout, we shall assume that $S^1 = \mathbb{R}/(2\pi\mathbb{Z})$, $B(S^1)$ is the set of all bounded, (Borel) measurable functions from $S^1$ to $\mathbb{R}$, $\|\cdot\|_\infty$ is the corresponding sup-norm and $\mathcal{M}(S^1)$ is the set of all (Borel) probability measures on $S^1$. In a slight abuse of notation, we shall also use $\|\cdot\|_\infty$ to denote the induced operator norm on $\mathcal{L}(S^1) =$ the set of all bounded, linear functions $T: B(S^1) \to B(S^1)$.

Suppose that $f$ is a smooth map on $S^1$ and define the deterministic system

$$(2) \qquad\qquad x_{n+1} = f(x_n).$$

(The smoothness assumptions are stronger than necessary—see Remark 1.) We are interested in the dynamics of the perturbed system

$$(3) \qquad\qquad X_{n+1}^\varepsilon = f(X_n^\varepsilon) + \varepsilon\sigma(X_n^\varepsilon)\chi_n \bmod 2\pi,$$

where $\chi_n$ is a family of i.i.d. standard normal random variables and $\sigma \in C^\infty(S^1)$. We assume there exist positive constants $\sigma_{lb}, \sigma_{ub}$ so that $\sigma_{lb} <$ $\sigma(x) < \sigma_{ub}$, $\forall x \in S^1$ and write $\mathbb{P}^x$ for the probability law of (3) given that $X_0^\varepsilon = x$. It is easy to see that $X_n^\varepsilon$ forms a (time homogeneous) Markov Chain on $S^1$ with transition operator $T^\varepsilon : B(S^1) \to B(S^1)$ given by

$$(4) \qquad T^\varepsilon \phi(x) = \mathbb{E}^x[\phi(X_1^\varepsilon)] = \mathbb{E}[\phi(f(x) + \varepsilon \sigma(x)\chi)] = \int_{S^1} \phi(y) \tilde{p}^\varepsilon(x, y)\, dy$$

for any $\phi \in B(S^1)$ and $x \in S^1$ where

$$\tilde{p}^\varepsilon(x, y) := \sum_{n \in \mathbb{Z}} p^\varepsilon(x, y + 2\pi n),$$

$$p^\varepsilon(x, y) := \frac{1}{\sqrt{2\pi}\varepsilon\sigma(x)} e^{-(y - f(x))^2/(2\sigma^2(x)\varepsilon^2)}.$$

Since $\tilde{p}^\varepsilon$ is smooth in both variables and $S^1$ is compact, $T^\varepsilon$ is a compact operator on $B(S^1)$ for any $\varepsilon > 0$ and hence, its spectrum, which we denote by $\sigma(T^\varepsilon)$, consists of a countable number of eigenvalues with 0 as the only possible limit point. The fact that $\|T^\varepsilon\|_\infty = 1$ of course implies that $\sigma(T^\varepsilon) \subset \{\lambda \in \mathbb{C} : |\lambda| \leq 1\}$. Moreover, $\inf\{\tilde{p}^\varepsilon(x, y) : x, y \in S^1\} > 0$ so that $X_n^\varepsilon$ has a unique, stationary distribution $\mu^\varepsilon$ and for any $x \in S^1$, $P^x(X_n \in \cdot)$ converges to $\mu^\varepsilon(\cdot)$ in total variation (see for instance, [4], Section 5.6). Therefore, $T^\varepsilon$ always has a simple eigenvalue at 1 and all other eigenvalues are strictly less than 1 in modulus. Our first result gives us asymptotic expressions for lower order eigenvalues. In what follows, $f^p$ denotes the $p$th iterate of $f$.

**Theorem 1.** *Suppose $f$ has a finite number of stable periodic orbits $P_i$, $i = 1, 2, \ldots, m_s$ and unstable periodic orbits $Q_i$ of period $q_i$, $i = 1, 2, \ldots, m_u$. Let $c_{s,i} = (f^{p_i})'(x_i)$ for some $x_i \in P_i$ and $c_{u,i} = (f^{q_i})'(y_i)$ for some $y_i \in Q_i$. Assume in addition that*

$$\lim_{n \to \infty} f^n(x) \in \bigcup_{i=1}^{m_s} P_i$$

*for all $x \in S^1 \setminus (\bigcup_{i=1}^{m_u} Q_i)$. Then for all $r > 0$, we can decompose $T^\varepsilon = T_{up}^\varepsilon + T_{lp}^\varepsilon$ so that for small $\varepsilon > 0$, we have $\|T_{lp}^\varepsilon\|_\infty < r$ and any eigenvalue of $T_{up}^\varepsilon$ with modulus greater than $r$ is of the form $\lambda + O(\varepsilon)$ with:*

(i) $\lambda = (c_{s,i}^j)^{1/p_i}$ *for some $i = 1, 2, \ldots, m_s$ and $j \geq 0$*

*or*

(ii) $\lambda = (|c_{u,i}|^{-1} c_{u,i}^{-j})^{1/q_i}$ *for some $i = 1, 2, \ldots, m_u$ and $j \geq 0$.*

Note that we include all branches of the $p_i$th and $q_i$th root in (i) and (ii).

Theorem 1 is really a statement about the limiting pseudoeigenvalues of $T^\varepsilon$. Recall that $\lambda$ is a $r$-pseudoeigenvalue of a compact operator $T$ if $\lambda \in \sigma(T + E)$ for some bounded linear operator $E$ with $\|E\| < r$ ([12], page 31). We also note that for any compact $T$, $\sigma(T) = \bigcap_{r>0} \sigma_r(T)$ where $\sigma_r(T)$ is the set of all $r$-pseudoeigenvalues of $T$ ([12], Theorem 4.3), but in this paper, we do not address the issue of taking double limits as $\varepsilon \to 0$ and $r \to 0$. Instead, for the remainder of this paper, we will say that $\lambda$ is a *limiting eigenvalue* of the operator $T^\varepsilon$ if $\forall r > 0$, $T^\varepsilon$ has a sequence of $r$-pseudoeigenvalues which converge to $\lambda$ as $\varepsilon \to 0$. Therefore, Theorem 1 states that $T^\varepsilon$ has limiting eigenvalues given by (i) and (ii) above.

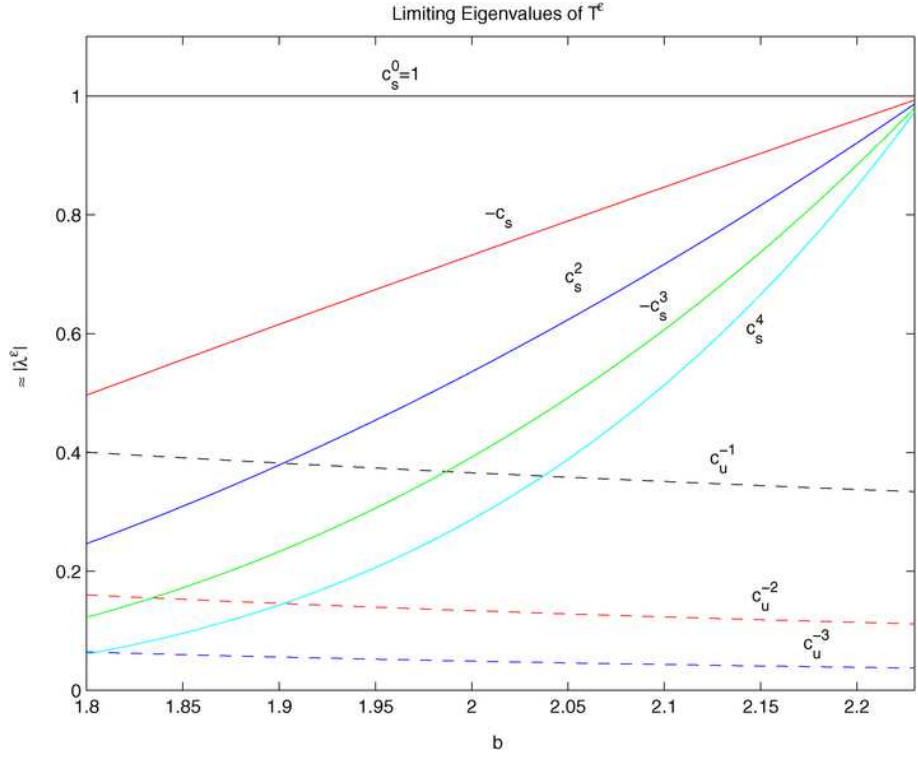

Fig. 1. *Modulus of the top limiting eigenvalues of $T^\varepsilon$ as given by Theorem 1 plotted against $b$ up to the first period doubling bifurcation point at $b = \sqrt{5} \approx 2.23$. The solid lines are powers of $c_s$ while the dashed lines are negative powers of $c_u$.*

We illustrate the results of Theorem 1 with a concrete example by taking $f(x) = x + 1 - b \sin x$ in the well-studied family of sine-circle maps (see, e.g., [5]). If $b > 1$, then $f$ has two fixed points $x_u, x_s \in (-\pi, \pi)$ with $c_u = f'(x_u) > 1$ for all $b > 1$ and $c_s = f'(x_s) \in (-1, 1)$ if and only if $b < b_c := \sqrt{5} \approx 2.23$. $f$ has no other periodic orbits for $b < b_c$. Therefore, if $1 < b < b_c$, Theorem 1 tells us that $T^\varepsilon$ has limiting eigenvalues $c_s^n$ and $c_u^{-(n+1)}$ for $n \geq 0$ (see Figure 1).

When $b = b_c$, $f'(x_s) = -1$ so that (2) undergoes a period doubling bifurcation with the appearance of a stable period two orbit. Figure 2 shows con-tributions to the spectrum of $T^\varepsilon$ coming from its two unstable fixed points $x_s$ and $x_u$ and stable period two orbit $P = \{x_1, x_2\}$. The contributions from the fixed point are of the form $|c_s|^{-1} c_s^{-n}$ or $|c_u|^{-1} c_u^{-n}$ where $c_{s,u} = f'(x_{s,u})$ and the contributions from $P$ are of the form $\sqrt{c^n}$ where $c = f'(x_1) f'(x_2)$ and we take both branches of the square root. This leads to the appearance of pairs of equal modulus eigenvalues in the bottom half of Figure 2. The quantitative change in the limiting eigenvalues of $T^\varepsilon$ near the deterministic bifurcation point $b_c$ motivates the following definition.

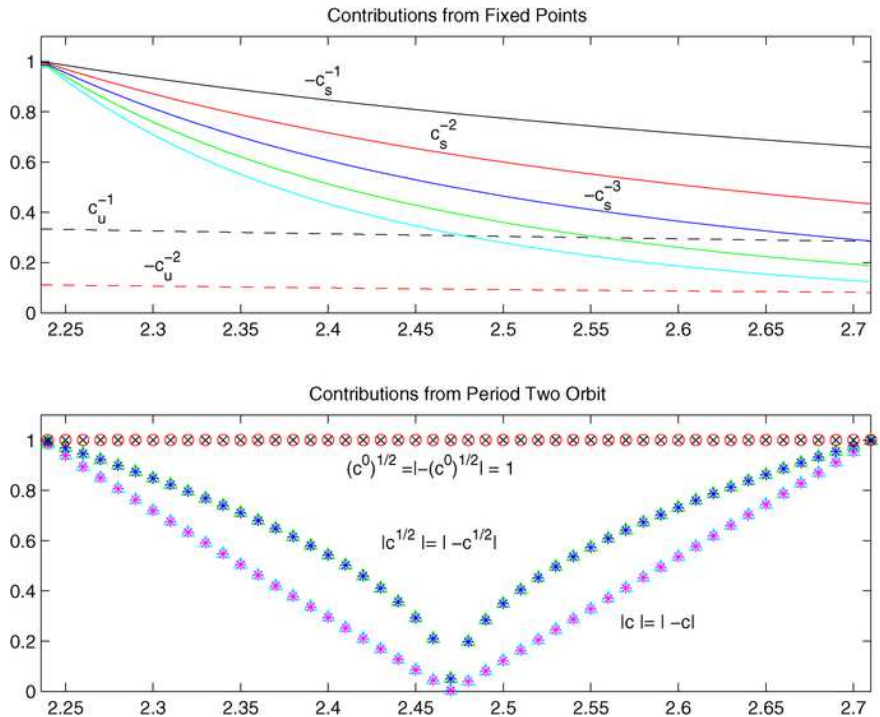

Fig. 2. *Modulus of the top limiting eigenvalues of $T^\varepsilon$ as $\varepsilon \to 0$ as given by Theorem 1 plotted against $b$ between the first period doubling bifurcation point at $b = \sqrt{5} \approx 2.23$ and the second near $b \approx 2.71$.*

DEFINITION 1. We call any change to the number of limiting eigenvalues of $T^\varepsilon$ with modulus 1 (as $b$ is varied) a $\lambda$-bifurcation.

Therefore, a $\lambda$-bifurcation occurs at $b_c$ with the appearance of a limiting eigenvalue at $-1$. As $b \nearrow 2.71$, $c \searrow -1$ and a second period doubling occurs in the deterministic system with the appearance of a stable period four orbit. Since a stable period four orbit yields four limiting eigenvalues which approach the unit circle as $\varepsilon \to 0$, another $\lambda$-bifurcation will occur near $b \approx 2.71$ as well. Figure 3 illustrates the complete $\lambda$-bifurcation scenario up to this second period doubling point.

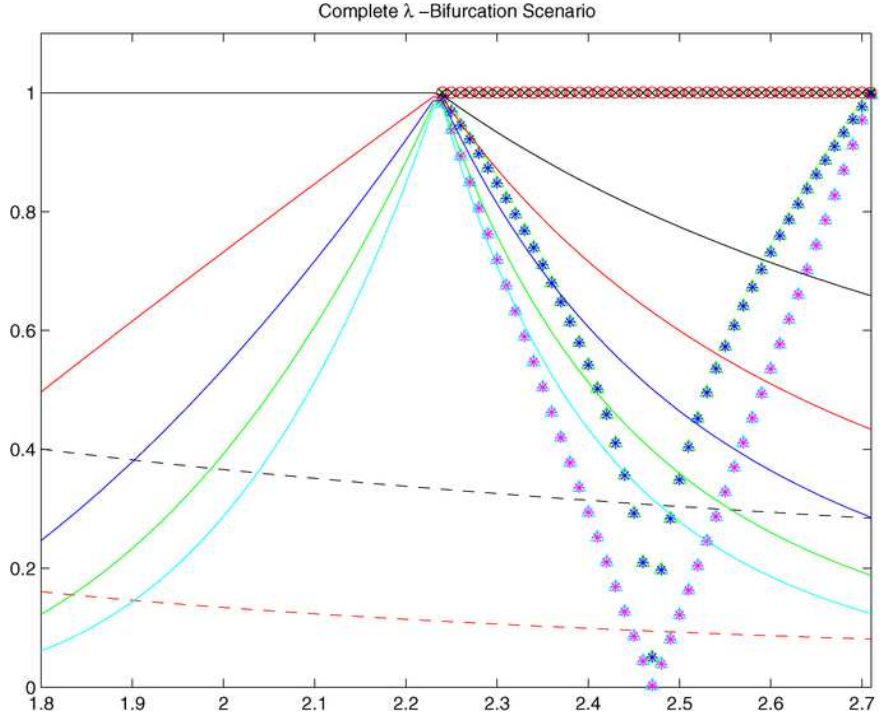

FIG. 3.   *Illustration of $\lambda$-bifurcation scenario for $f(x) = x + 1 - b\sin x$, $1.8 < b < 2.71$. See also Figures 1 and 2.*

Our second result deals with eigenvectors. We use the notation $H_n$ to denote the $n$th Hermite polynomial (see [8] for definitions) and $h_n(x) = e^{-x^2}H_n(x)$ for the corresponding Hermite function.

THEOREM 2.   *If $P = \{x_1, \ldots, x_p\}$ is a stable periodic orbit of $f$ with $c_s = (f^p)'(x_1)$, then the eigendensities of $T^\varepsilon$ corresponding to the limiting eigenvalue $(c_s^n)^{1/p}$ are of the form*

$$\sum_{j=1}^{p} a_j \phi_{s,n,j}(x) + O(\varepsilon)$$

*for some constants $a_j$ where $\phi_{s,n,j}(x) = h_n(\alpha_j(x - x_j)/\varepsilon)$ with $\alpha_j$ an explicit constant depending on $c_i = f'(x_i)$, $i \neq j$ and $\sigma(x_i)$, $i = 1, 2, \ldots, p$. If $Q = \{y_1, \ldots, y_q\}$ is an unstable periodic orbit with $c_u = (f^q)'(y_1)$, then the eigenfunctions corresponding to the limiting eigenvalue $(|c_u|^{-1}c_u^{-n})^{1/q}$ are of*

*the form*

$$\sum_{j=1}^{q} b_j \phi_{u,n,j}(x) + O(\varepsilon)$$

*for some constants $b_j$ where $\phi_{u,n,j}(x) = h_n(\beta_j(x - y_j)/\varepsilon)$ and $\beta_j$ is an explicit constant depending on $c_i = f'(x_i)$, $i \neq j$, and $\sigma(x_i)$ for $i = 1, 2, \ldots, p$.*

See Theorem 3 in Section 3 and Theorem 8 in Section 6 for more details in the period one and two case, including formulas for $\alpha_j$ and $\beta_j$. Figure 4 illustrates the limiting invariant densities ($n = 0$ eigendensities) from Theorem 2 for two different parameter values in the sine-circle example discussed above. If we take $b$ just past the second period doubling bifurcation point (when the deterministic system has stable period 4 behavior), the amount of "humps" in the limiting invariant density for the perturbed system will again double. Therefore, we can see that the shape of the invariant density is greatly affected by the number of limiting eigenvalues near the unit circle. Qualitative changes to the shape of invariant densities are often called $P$-bifurcations (see [2] for examples). Further connections between $\lambda$-bifurcations and $P$-bifurcations may also be an interesting question for future work.

The reason we are only able to calculate eigendensities for one set of eigenvalues and eigenfunctions for the other is a direct consequence of the structure of $T^\varepsilon$. We now move on to discuss this structure and give a detailed proof of our results in the stable period one case, returning to the general case in Section 6.

## 3. One stable and one unstable fixed point: basic setup and main results.

Throughout this section, we assume that $f$ has two fixed points $x_s, x_u$ satisfying $c_s := f'(x_s) \in (-1, 1)$ and $c_u := f'(x_u) \notin [-1, 1]$ with the property that $f^n(x) \to x_s$ as $n \to \infty$, $\forall x \in S^1$, $x \neq x_u$. This corresponds to studying perturbations of (2) in a regime of stable period one behavior. We set $\sigma_s = \sigma(x_s)$, $\sigma_u = \sigma(x_u)$, $\alpha = \sqrt{(1 - c_s^2)/(2\sigma_s^2)}$ and $\beta = \sqrt{(c_u^2 - 1)/(2\sigma_u^2)}$. $H_n$ and $h_n$ are is in Theorem 2. In the language of Section 2, the following result says that $T^\varepsilon$ has limiting eigenvalues $c_s^j$ and $|c_u|^{-1} c_u^{-j}$.

THEOREM 3. *Suppose that $f$ is a smooth map on $S^1$ with stable fixed point $x_s$ and unstable fixed point $x_u$. In addition, assume that $f^n(x) \to x_s$ for all $x \in S^1 \setminus \{x_u\}$. Then for any $r > 0$, $\exists \varepsilon_r, L_r, K_r > 0$, so that $\forall \varepsilon < \varepsilon_r$, we can write $T^\varepsilon = T_{up}^\varepsilon + T_{lp}^\varepsilon$ where*

$$\|T_{lp}^\varepsilon\|_\infty < r$$

*and any $\lambda \in \sigma(T_{up}^\varepsilon)$ with $|\lambda| > r$ is a simple eigenvalue of one of the two forms:*

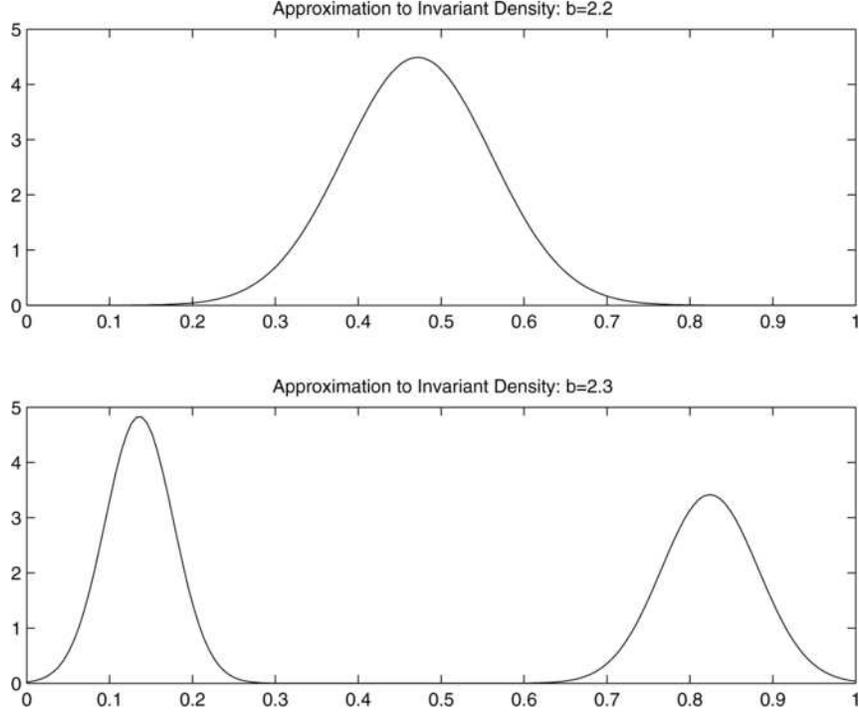

Fig. 4. *Approximations from Theorem 2 for the limiting invariant densities of (3) in the case when $f(x) = x + 1 - b \sin x$. In the top, we take $b = 2.2$ so that $f$ has stable fixed point $x_s \approx 0.47$ and in the bottom, we take $b = 2.3$ so that $f$ has stable period two orbit $P \approx \{0.14, 0.82\}$.*

    (i)  $\lambda = c_s^j + \lambda_{s,j,1}^\varepsilon \varepsilon$

*or*

    (ii)  $\lambda = |c_u|^{-1} c_u^{-j} + \lambda_{u,j,1}^\varepsilon \varepsilon$

*for some $j \geq 0$ with $c_s = f'(x_s)$, $c_u = f'(x_u)$, and $\max(|\lambda_{s,j,1}^\varepsilon|, |\lambda_{u,j,1}^\varepsilon|) \leq L_r$. All eigendensities corresponding to $\lambda$ as in (i) are multiples of*

$$(5) \qquad \left[ h_j\left( \frac{\alpha(x - x_s)}{\varepsilon} \right) + \varepsilon \psi_{s,j}^\varepsilon \left( \frac{x - x_s}{\varepsilon} \right) \right] \mathbf{1}_{V_3}(x),$$

*where $V_3$ is a neighborhood of $x_s$ and $\psi_{s,j}^\varepsilon$ has the property that*

$$\sup_{x \in \mathbb{R}} (|\psi_{s,j}^\varepsilon(x)| e^{kx^2}) < K_r$$

*for some $k > 0$. All eigenfunctions corresponding to $\lambda$ as in (ii) are multiples of*

$$(6) \qquad \left[ h_j\left( \frac{\beta(x - x_u)}{\varepsilon} \right) + \varepsilon \psi_{u,j}^\varepsilon \left( \frac{x - x_u}{\varepsilon} \right) \right] \mathbf{1}_{V_1}(x)$$

*with $V_1$ a neighborhood of $x_u$ and*

$$\sup_{x\in\mathbb{R}}(|\psi_{u,j}^{\varepsilon}(x)|e^{kx^2}) < K_r.$$

REMARK 1. We assume $f$ is smooth only for convenience. In fact, if as long as $f \in C^{n+1}$ for some $n \geq 1$ in a neighborhood of each fixed point (and possibly discontinuous elsewhere), it can be shown that $\exists \lambda_{s,j,k} \in \mathbb{C}$, $k = 1, 2, \ldots, n-1$ such that the eigenvalues in (i) have asymptotic expansions

$$\lambda = c_{s,j}^j + \sum_{k=2}^{n-1} \varepsilon^k \lambda_{s,j,k} + \lambda_{s,j,n}^{\varepsilon} \varepsilon^n$$

with $|\lambda_{s,j,n}^{\varepsilon}| \leq L_r$, $\forall \varepsilon < \varepsilon_r$ and similarly for the eigenvalues in (ii). We leave the details to the reader. The key is in getting higher order terms out of Lemma 6. See also Remark 5.

REMARK 2. The coefficients $\lambda_{s,j,k}$ in Remark 1 can be calculated if $\sigma$ is constant using standard properties of Hermite polynomials (see [8] for a list of some of these properties). In particular, it can be shown that in this case, $\lambda_{s,j,1} = 0$ so that convergence is in fact of order $\varepsilon^2$.

The starting point for the proof of Theorem 3 is Proposition 1 below which gives us a way of splitting up the circle into regions determined by the different actions of $f$. In what follows, $d$ denotes the standard quotient metric on $S^1$ induced by the Euclidean metric on $\mathbb{R}$ and $B_\delta(x) = \{y \in S^1 : d(x, y) < \delta\}$, $\delta > 0$.

PROPOSITION 1. *There exist neighborhoods $V_1 := B_{\delta_u}(x_u)$, $V_3 := B_{\delta_s}(x_s)$, and constants $\eta > 0$, $N \in \mathbb{N}$ such that:*

(i) $d(f(x), V_1) > \eta$ *for every* $x \notin V_1$.
(ii) $d(f(x), V_3^c) > \eta$ *for every* $x \in V_3$.
(iii) *For every* $x \in V_2 := S^1 \setminus (V_1 \cup V_3)$, *we have* $f^n(x) \in V_3$, $\forall n \geq N$.

PROOF. To prove (i), we first choose $\delta' > 0$ so that $d(f(x), x_u) > \gamma_u d(x, x_u)$, $\forall x \in B_{\delta'}(x_u)$ with $\gamma_u > 1$ and let $K = S^1 \setminus B_{\delta'}(x_u)$. $K$ is compact so that $f(K)$ is compact and therefore, we can find a $\delta_u \in (0, \delta')$ so that $f(K) \subset S^1 \setminus B_{2\delta_u}(x_u)$. Then if $x \notin B_{\delta_u}(x_u) =: V_1$ either: $x \in K$, in which case $d(f(x), x_u) \geq 2\delta_u$ or $x \in B_{\delta'} \setminus B_{\delta_u}(x_u)$, in which case $d(f(x), x_u) > \gamma_u d(x, x_u) > \gamma_u \delta_u$. This implies that $d(f(x), V_1) > \eta_1$ for every $x \notin V_1$ with $\eta_1 := \min(\gamma_u - 1, 1)\delta_u$. (ii) follows from a similar argument and the fact that $|f'(x_s)| < 1$ while (iii) follows directly from the assumption that $f^n(x) \to x_s$ as $n \to \infty$, $\forall x \in V_2$, and the compactness of $V_2$. $\square$

Motivated by the above proposition, we write $\phi \in B(S^1)$ as $\phi = \phi_1 + \phi_2 + \phi_3$ where $\phi_i = \phi \mathbf{1}_{V_i}$, $i = 1, 2, 3$, and then decompose $T^\varepsilon$ into operators $T_{ij}^\varepsilon : B(V_j) \to B(V_i)$ defined by

$$(7) \qquad T_{ij}^\varepsilon \phi(x) = \mathbf{1}_{V_i}(x) \mathbb{E}[(\phi \mathbf{1}_{V_j})(f(x) + \varepsilon \sigma(x) \chi)] = \int \phi(y) \tilde{p}_{ij}^\varepsilon(x, y) \, dy,$$

where

$$\tilde{p}_{ij}^\varepsilon(x, y) := \mathbf{1}_{V_i}(x) \tilde{p}^\varepsilon(x, y) \mathbf{1}_{V_j}(y).$$

We can think of $T^\varepsilon$ acting on $\phi \in B(S^1)$ via matrix multiplication:

$$T^\varepsilon \phi(x) \Leftrightarrow \begin{bmatrix} T_{11}^\varepsilon & T_{12}^\varepsilon & T_{13}^\varepsilon \\ T_{21}^\varepsilon & T_{22}^\varepsilon & T_{23}^\varepsilon \\ T_{31}^\varepsilon & T_{32}^\varepsilon & T_{33}^\varepsilon \end{bmatrix} \begin{bmatrix} \phi_1 \\ \phi_2 \\ \phi_3 \end{bmatrix}(x).$$

(This can be made precise by use of inclusion/restriction operators.) Recalling that $T^\varepsilon$ is the transition operator for (3), (7) implies that we can informally think of $T_{ij}^\varepsilon$ as providing information about movement from $V_i$ to $V_j$.

If we take $\varepsilon = 0$, we obtain the deterministic system (2). Furthermore, Proposition 1 implies that the transition operator $T^0 \phi(x) = \phi(f(x))$ has the "upper triangular" decomposition

$$T^0 = \begin{bmatrix} T_{11}^0 & T_{12}^0 & T_{13}^0 \\ 0 & T_{22}^0 & T_{23}^0 \\ 0 & 0 & T_{33}^0 \end{bmatrix}$$

with the additional property that $(T_{22}^0)^n = 0$, $\forall n \geq N$. With noise in the system, we cannot hope for such good fortune as there is always a small probability of movement between regions. We can, however, obtain bounds on the probabilities of such events, as the next three lemmas illustrate. At this point, we also introduce the notation $B(U, V) =$ set of all bounded measurable functions from $U \subset S^1$ to $V \subset S^1$ and $B(U) = B(U, U)$. Let $\| \cdot \|_{\infty, U, V}$ denote the corresponding sup-norm. We will simply write $\|\phi\|_\infty$ when the domain and range of $\phi$ are clear. We shall use the same notation and caveats when referring to the induced operator norm on $\mathcal{L}(B(U), B(V))$ = set of all bounded linear, operators from $B(U)$ to $B(V)$ [with $\mathcal{L}(B(U)) = \mathcal{L}(B(U), B(U))$]. For instance, In the following lemmas, we have $T_{ij}^\varepsilon : B(V_j) \to B(V_i)$ so we write $\|T_{ij}^\varepsilon\|_\infty$ for $\|T_{ij}^\varepsilon\|_{\infty, V_j, V_i}$.

LEMMA 1. *There exist constants $M, K > 0$ so that*

$$\|T_{ij}^\varepsilon\|_\infty \leq M \varepsilon e^{-K/\varepsilon^2}$$

*for every $i > j$.*

PROOF. Clearly, $\|T_{ij}^{\varepsilon}\|_{\infty} = \sup_{x \in V_i} \mathbb{P}(f(x) + \varepsilon \sigma(x) \chi \in V_j)$. If $x \in V_3$ and $j = 1, 2$, then (ii) in Proposition 1 implies that

$$\mathbb{P}(f(x) + \varepsilon \sigma(x) \chi \in V_j) \leq \mathbb{P}(d(f(x) + \varepsilon \sigma(x) \chi, f(x)) > \eta).$$

Similarly, if $x \in V_2$ and $j = 1$, (i) in Proposition 1 implies that

$$P(f(x) + \varepsilon \sigma(x) \chi \in V_1) \leq P(d(f(x) + \varepsilon \sigma(x) \chi, f(x)) > \eta).$$

The result is then a direct consequence of the next lemma. □

LEMMA 2. *For any $a, \varepsilon > 0$ and $x \in S^1$,*

$$\mathbb{P}(d(f(x) + \varepsilon \sigma(x) \chi, f(x)) > a) \leq \sqrt{\frac{2\sigma_{ub}}{\pi a^2}} \varepsilon e^{-a^2/(2\sigma_{ub}^2 \varepsilon^2)}.$$

PROOF. Follows from standard normal distribution tail estimates. □

LEMMA 3. *There exist positive constants $M_N, K_N$ such that*

$$\|(T_{22}^{\varepsilon})^{N+1}\|_{\infty} \leq M_N \varepsilon e^{-K_N/\varepsilon^2},$$

*where $N$ is the same constant as in* (iii) *of Proposition 1.*

PROOF. From (iii) in Proposition 1 we know that $f^N(x) \in V_3$, $\forall x \in V_2$ so that $d(f^{N+1}(x), V_2) > \eta$ by (ii). Therefore,

$$\|(T_{ij}^{\varepsilon})^{N+1}\|_{\infty} \leq \sup_{x \in V_2} \mathbb{P}^x(d(X_{N+1}^{\varepsilon}, f^{N+1}(x)) > \eta).$$

Since

$$d(X_{N+1}^{\varepsilon}, f^{N+1}(x)) \leq \sum_{i=0}^{N} L^{N-i} d(X_{i+1}^{\varepsilon}, f^i(x)),$$

where $L = \sup |f'(x)|$, the result then follows from independence and Lemma 2. □

With these results in hand, we are ready to give the following:

PROOF OF THEOREM 3. For any $\varepsilon > 0$, we can write $T^{\varepsilon} = T_{up}^{\varepsilon} + T_{lp}^{\varepsilon}$ where

$$T_{up}^{\varepsilon} := \begin{bmatrix} T_{11}^{\varepsilon} & T_{12}^{\varepsilon} & T_{13}^{\varepsilon} \\ 0 & T_{22}^{\varepsilon} & T_{23}^{\varepsilon} \\ 0 & 0 & T_{33}^{\varepsilon} \end{bmatrix}$$

and

$$T_{lp}^{\varepsilon} := \begin{bmatrix} 0 & 0 & 0 \\ T_{21}^{\varepsilon} & 0 & 0 \\ T_{31}^{\varepsilon} & T_{32}^{\varepsilon} & 0 \end{bmatrix}.$$

If $\varepsilon$ is sufficiently small, then Lemma 1 implies that $\|T^{\varepsilon}_{lp}\|_{\infty} \leq M\varepsilon e^{-K/\varepsilon^2} < r$. Since $T^{\varepsilon}_{up}$ is upper triangular, its spectrum is included in the union of the spectra of the diagonal operators $T^{\varepsilon}_{ii}$, $i = 1, 2, 3$. But by Lemma 3 and the fact that $\|T^{\varepsilon}_{22}\|_{\infty} \leq 1$, the spectral radius of $T^{\varepsilon}_{22}$ can be made less than $r$ by shrinking $\varepsilon$ if necessary so that any eigenvalue of $T^{\varepsilon}$ with modulus greater than $r$ must be in $\sigma(T^{\varepsilon}_{11})$ or $\sigma(T^{\varepsilon}_{33})$. Furthermore, because of the upper triangular structure, we know that if $\phi_3$ is an eigendensity of $T^{\varepsilon}_{33}$, then

$$\phi(x) := \begin{cases} \phi_3(x), & x \in V_3, \\ 0, & x \notin V_3, \end{cases}$$

is an eigendensity of $T^{\varepsilon}_{up}$ and similarly, the eigenfunctions of $T^{\varepsilon}_{11}$ yields eigenfunctions of $T^{\varepsilon}_{up}$. Therefore, the proof of Theorem 3 will be complete if we can show that all eigenvalues of $T^{\varepsilon}_{33}$ with modulus larger than $r$ are of the form (i) with corresponding eigendensities (5) and all eigenvalues of $T^{\varepsilon}_{11}$ with modulus larger than $r$ are of the form (ii) with corresponding eigenfunctions (6). Calculating the spectra of these operators turns out to be a difficult task and is of interest in its own right. We therefore, dedicate the next two sections to this analysis and note that Theorems 4 and 6 in Sections 4 and 5, respectively, give the results necessary for the completion of this proof.  □

## 4. The local story near a stable fixed point.

The essential conclusions from our work in this section are contained in the following theorem which provides us with the necessary information we need about the action of $T^{\varepsilon}$ near a stable fixed point of $f$.

THEOREM 4.   *For any $r > 0$, $\exists \varepsilon_{s,r}, L_{s,r} > 0$ so that $\forall \varepsilon < \varepsilon_{s,r}$, any eigenvalue of $T^{\varepsilon}_{33}$ in $B(V_3)$ is a simple eigenvalue of the form*

$$\lambda^{\varepsilon}_{s,j} = c^j_s + \varepsilon \lambda^{\varepsilon}_{s,j,1}$$

*for some $j \geq 0$ with $|\lambda^{\varepsilon}_{s,j,1}| \leq L_{s,r}$, $\forall j \geq 0, \varepsilon < \varepsilon_{s,r}$. Furthermore, $\exists K_{s,r} > 0$ such that the eigendensities of $T^{\varepsilon}_{33}$ corresponding to $\lambda^{\varepsilon}_{s,j}$ are multiples of*

$$\left[ h_j\left( \frac{\alpha(x - x_s)}{\varepsilon} \right) + \varepsilon \psi^{\varepsilon}_{s,j}\left( \frac{x - x_s}{\varepsilon} \right) \right] \mathbf{1}_{V_3}(x)$$

*with $h_j, \alpha$ as in Theorem 3 and*

$$\sup_{x \in \mathbb{R}} (|\psi^{\varepsilon}_{s,j}(x)| e^{kx^2}) \leq K_{s,r}$$

*for all $\varepsilon < \varepsilon_{s,r}$ and some $k > 0$.*

Before delving into the details of the proof [which are rather complicated due to the singular nature of the perturbation in (3)], we first provide some motivation. We identify $S^1$ with $[-\pi/2, \pi/2]$ and $x_s$ with $0$ so that $V_3 =$ $(-\delta_s, \delta_s)$, $f(0) = 0$, $c_s = f'(0)$, and $T^\varepsilon_{33} : B((-\delta_s, \delta_s)) \to B((-\delta_s, \delta_s))$. As we will be working on $\mathbb{R}$ for the remainder of this section, we will now set $B_\delta(x) := \{ y \in \mathbb{R} : |x - y| < \delta \}$.

We can informally think of (3) as a Markov Chain on $\mathbb{R}$ [with $f(0) = 0$] and re-scale space near the origin. In other words, we look at the chain $Y^\varepsilon_n := X^\varepsilon_n / \varepsilon$. Then

$$Y^\varepsilon_{n+1} = X^\varepsilon_{n+1} / \varepsilon = \varepsilon^{-1} f(X^\varepsilon_n) + \sigma(X^\varepsilon_n) \chi_n = f^\varepsilon(Y^\varepsilon_n) + \sigma^\varepsilon(Y^\varepsilon_n) \chi_n,$$

where $f^\varepsilon(x) = \varepsilon^{-1} f(\varepsilon x)$ and $\sigma^\varepsilon(x) = \sigma(\varepsilon x)$. Since $f^\varepsilon(x) \to c_s x$ and $\sigma^\varepsilon(x) \to \sigma(0)$ as $\varepsilon \to 0$, we expect that the dynamics of $Y^\varepsilon_n$ should be closely approximated by the dynamics of

$$Y_{n+1} = c_s Y_n + \sigma(0) \chi_n$$

for small values of $\varepsilon$.

This limit is nondegenerate and describes a simple autoregressive scheme. It can easily be verified that $Y_n$ has an invariant measure $\mu$ defined by $\mu(dx) = \rho_s(x) \, dx$ where

$$\rho_s(x) := \frac{\alpha}{\sqrt{\pi}} e^{-(\alpha y)^2}$$

with $\alpha = \sqrt{(1 - c_s^2)/(2\sigma_0^2)}$. Therefore, the space $L^2(\mu)$ provides a natural setting for investigating the spectrum of the transition operator, $T_s$, for $Y_n$ defined by

$$T_s \phi(x) = \int \phi(y) p_s(x, y) \, dy \qquad \text{with } p_s(x, y) = \frac{1}{\sqrt{2\pi} \sigma(0)} e^{-(y - c_s x)^2 / (2\sigma^2(0))}.$$

For the following results, we write $\| \cdot \|_2$ for $\| \cdot \|_{L^2(\mu)}$.

LEMMA 4. *$T_s$ acts as a bounded, self-adjoint operator on $L^2(\mu)$ with $\|T_s\|_2 = 1$.*

PROOF. Apply the Cauchy–Schwarz inequality and note that

(8) $$\rho_s(x) p_s(x, y) = \rho_s(y) p_s(y, x)$$

for all $x, y \in \mathbb{R}$. □

Since $T_s$ is a self-adjoint operator on $L^2(\mu)$, we know that we can find a complete, orthonormal set (CONS) of eigenfunctions for $T_s$ in $L^2(\mu)$. The following lemma identifies these functions.

LEMMA 5. *The eigenvalues of $T_s$ in $L^2(\mu)$ are given by $c_s^n$, $n \geq 0$ and the corresponding eigenfunctions are multiples of $\phi_{s,n}(x) = H_n(\alpha x)$ where $H_n$ is the $n$th Hermite polynomial.*

PROOF. Using the generating function definition of Hermite polynomials as the functions satisfying

$$\sum_{n=0}^{\infty} \frac{H_n(x)z^n}{n!} = e^{-z^2 + 2xz}$$

and the fact that $\mathbb{E}(e^{t\chi}) = e^{t^2/2}$ for any $t \in \mathbb{R}$ when $\chi$ is standard normal, we obtain

$$(9) \qquad T_s\left(\sum_{n=0}^{\infty} \frac{H_n(\alpha x)}{n!} z^n\right) = \sum_{n=0}^{\infty} \frac{H_n(\alpha x)}{n!}(c_s z)^n.$$

Since

$$\int H_n(\alpha x) H_m(\alpha x)\, d\mu(x) = \delta_{n,m} 2^n n!,$$

the partial sums $S_N(x) := \sum_{n=0}^{N} \frac{H_n(\alpha x)z^n}{n!}$ form a Cauchy sequence in $L^2(\mu)$ and hence, (9) and the continuity of $T_s$ imply that the $c_s^n$ are in fact eigenvalues with corresponding eigenfunctions $\phi_{s,n}$. Since the $H_n$ form a CONS in $L^2(\nu)$ where $\nu$ is defined by $\nu(dx) = \sqrt{\frac{2}{\pi}} e^{-x^2}\, dx$ (see [8]), the $\phi_n$ form a CONS in $L^2(\mu)$, which proves the result. $\square$

REMARK 3. For any (Borel) measure $m$, define the measure $mT_s(A) = \int T_s \mathbf{1}_A(x) m(dx) = \int \int_A p_s(x,y)\, dy\, m(dx)$ for all (Borel) measurable sets $A$. If we let $\phi_{s,n}^*(dx) = \phi_{s,n}(x)\mu(dx)$, then by (8) and Lemma 5, a quick calculation shows that $\forall A \in \mathcal{B}$, $\phi_{s,n}^* T_s(A) = c^n \phi_{s,n}^*$. Therefore, $T_s$ has eigenmeasures $\phi_{s,n}^*$ and eigendensities $\phi_{s,n}(x)\rho_s(x) = h_n(\alpha x)$.

Lemma 5 explains the limits in Theorem 4. The next two subsections give the technical arguments.

### 4.1. Expansion of the transition operator.

Instead of directly extending $T_{33}^\varepsilon$ to an operator on $B(\mathbb{R})$, we first define the family of weighted sup-norms:

$$(10) \qquad \|\phi\|_{A,k} = \sup_{x \in A} \frac{|\phi(x)|}{v_k(x)}$$

with $v_k(x) = e^{kx^2}$ and let $W_{A,k} = \{\phi \in \mathcal{B}(A) : \|\phi\|_{A,k} < \infty\}$ where $\mathcal{B}(A)$ denotes the set of all (Borel) measurable functions on $A$. It is easy to show that for any set $A$, $W_{A,k}$ along with the $\|\cdot\|_{A,k}$ norm is a Banach space [if $k = 0$, $W_{A,k} = B(A)$ and $\|\cdot\|_{A,k}$ is just the sup-norm]. When $A = \mathbb{R}$, we drop the $A$ dependence and write $\|\cdot\|_k$ and $W_k$ for $\|\cdot\|_{A,k}$ and $W_{A,k}$, respectively. Again in slight abuse of notation, we shall use $\|\cdot\|_k$ to refer to the operator norm on $\mathcal{L}(W_k) = $ set of all bounded linear operators on $W_k$ as well.

Since $V_3$ is bounded, the $\|\cdot\|_{V_3,k}$ norms are equivalent for all $k$ and hence, the spectrum of $T_{33}^\varepsilon$ will not depend on our use of norm. Therefore, we let $k > 0$ and use the norm $\|\cdot\|_{V_3,k/\varepsilon^2}$ on $B(V_3)$. Re-scaling can be done by applying the operator $U^\varepsilon : B(V_3^\varepsilon) \to B(V_3)$, $V_3^\varepsilon := (-\delta_s/\varepsilon, \delta_s/\varepsilon)$, defined by $U^\varepsilon \phi(x) := \phi(x/\varepsilon)$ and setting $T_3^\varepsilon := (U^\varepsilon)^{-1} \circ T_{33}^\varepsilon \circ U^\varepsilon$. Then

$$
\begin{aligned}
T_3^\varepsilon \phi(x) &= \int \phi(y/\varepsilon) \tilde{p}_{i,j}^\varepsilon(\varepsilon x, y)\, dy \\
&= \mathbf{1}_{V_3^\varepsilon}(x) \int_{V_3^\varepsilon} \phi(y) \varepsilon \tilde{p}^\varepsilon(\varepsilon x, \varepsilon y)\, dy
\end{aligned}
\tag{11}
$$

and the spectrum of $T_3^\varepsilon$ in $W_{V_3^\varepsilon, k}$ will be the same as the spectrum of $T_{33}^\varepsilon$ in $W_{V_3, k/\varepsilon^2}$. Finally, we can extend $T_3^\varepsilon$ to an operator on $\mathbb{R}$ via (11) and consider the spectrum of the resulting operator, $T_s^\varepsilon : W_k \to W_k$. Note that if $T_s^\varepsilon \phi(x) = \lambda \phi(x)$ for some $\phi \in W_k$ and $\lambda \neq 0$, then $\hat{\phi} := \phi|_{V_3^\varepsilon} \in W_{k, V_3^\varepsilon}$ satisfies $T_3^\varepsilon \hat{\phi}(x) = \lambda \hat{\phi}(x)$. Conversely, if $T_3^\varepsilon \hat{\phi}(x) = \lambda \hat{\phi}(x)$ for some $\hat{\phi} \in W_{k, V_3^\varepsilon}$ and $\lambda \neq 0$, then we can extend $\hat{\phi}$ to a function $\phi \in W_k$ such that $T_s^\varepsilon \phi(x) = \lambda \phi(x)$ by setting $\phi(x) = \hat{\phi}(x)$ for all $x \in V_3^\varepsilon$ and $\phi(x) = 0$, $\forall x \notin V_3^\varepsilon$. Therefore, the nonzero part of the spectrum of $T_3^\varepsilon$ will not be affected by this extension.

The $W_k$ spaces are large enough to include the eigenfunctions $\phi_{s,n}$ of $T_s$ and hence, are a good candidate space for studying the convergence of $T_s^\varepsilon$ to $T_s$. In fact, we can show that:

THEOREM 5. *There exists a $k_s > 0$ so that $T_s^\varepsilon = T_s + O(\varepsilon)$ in $\mathcal{L}(W_k)$ for all $k \in (0, k_s)$.*

The essential ideas in the proof are the expansion of the main part of $\varepsilon \tilde{p}^\varepsilon(\varepsilon \cdot, \varepsilon \cdot)$ about $p_s(\cdot, \cdot)$ and the use of the weight functions $v_k(\cdot)$ to control the growth of error terms. We begin with an expansion for the transition densities. Recall that

$$
T_s^\varepsilon \phi(x) = \int \phi(y) \mathbf{1}_{V_3}(\varepsilon x) \mathbf{1}_{V_3}(\varepsilon y) \varepsilon \tilde{p}^\varepsilon(\varepsilon x, \varepsilon y)\, dy,
$$

where

$$
\varepsilon \tilde{p}^\varepsilon(\varepsilon x, \varepsilon y) = \sum_{n \in \mathbb{Z}} \frac{1}{\sqrt{2\pi} \sigma(\varepsilon x)} e^{-(y + 2\pi n/\varepsilon - f^\varepsilon(x))^2/(2\sigma^2(\varepsilon x))}
$$

and $f^\varepsilon(x) = \varepsilon^{-1} f(\varepsilon x)$. We write

$$
p_m^\varepsilon(x, y) := \varepsilon p^\varepsilon(\varepsilon x, \varepsilon y) = \frac{1}{\sqrt{2\pi} \sigma(\varepsilon x)} e^{-(y - f^\varepsilon(x))^2/(2\sigma^2(\varepsilon x))}
$$

for the main part of the transition density for $T_s^\varepsilon$. Since $f^\varepsilon(x) \to c_s x$ and $\sigma(\varepsilon x) \to \sigma(0)$ as $\varepsilon \to 0$, we have

$$
p_m^\varepsilon(x, y) \to \frac{1}{\sqrt{2\pi} \sigma(0)} e^{-(y - c_s x)^2/(2\sigma^2(0))} = p_s(x, y)
$$

pointwise as $\varepsilon \to 0$. We also note that

$$(12) \qquad p_m^\varepsilon(x, y) \leq \frac{\sigma_{ub}}{\sigma_{lb}} \bar{p}_m^\varepsilon(x, y)$$

for all $x, y \in \mathbb{R}$ where

$$\bar{p}_m^\varepsilon(x, y) = \frac{1}{\sqrt{2\pi}\sigma_{ub}} e^{-(y - f^\varepsilon(x))^2/(2\sigma_{ub}^2)}.$$

LEMMA 6. *Let $\epsilon > 0$ and $U$ be a bounded subset of $\mathbb{R}$. Then $\forall x \in U/\varepsilon$ and $y \in \mathbb{R}$, we have*

$$p_m^\epsilon(x, y) = p_s(x, y) + \epsilon g_1(x, y) p_s(x, y) + \epsilon^2 R_\epsilon(x, y),$$

*where $g_1$ is a polynomial in $x, y$ while*

$$(13) \qquad |R_\epsilon(x, y)| \leq g_r^\epsilon(x, y)(\bar{p}_m^\epsilon(x, y) + p_s(x, y))$$

*with $g_r^\epsilon$ a polynomial in $\epsilon$, $y$, and $x$.*

PROOF. If we take a second order Taylor expansion of $g(z) = e^{-z}$ about $z = z_0$, we obtain

$$e^{-z} = e^{-z_0} - e^{-z_0}(z - z_0) + \frac{R(z, z_0)}{2}(z - z_0)^2$$

with $|R(z, z_0)| \leq \max\{e^{-z}, e^{-z_0}\} \leq e^{-z} + e^{-z_0}$, $\forall z \in \mathbb{R}$. Letting $z = \frac{(y - f^\epsilon(x))^2}{2\sigma^2(\varepsilon x)}$, $z_0 = \frac{(y - c_s x)^2}{2\sigma^2(0)}$, using the Taylor expansions of $f$ and $\sigma^{-1}$ about 0, and applying (12) then yields the result. $\square$

We will also need to establish some bounds on growth rates in $W_k$. To this end, define the linear operator

$$T_m^\varepsilon \phi(x) = \mathbf{1}_{V_3^\varepsilon}(x) \int_{V_3^\varepsilon} \phi(y) p_m^\varepsilon(x, y)\, dy$$

for $\phi \in W_k$. Let $I_\delta := (-\delta, \delta)$, $I_\delta^\varepsilon := (-\delta/\varepsilon, \delta/\varepsilon)$, and $M_\delta := \sup_{x \in I_\delta}\{|f'(x)|\}$.

LEMMA 7. *For any $k < \frac{1}{2\sigma_{ub}^2}$ and $j \in \mathbb{N}$, we have*

$$\int |y|^j v_k(x) p_s(x, y)\, dy \leq q(x) e^{l_k c_s^2 x^2}$$

*for all $x \in \mathbb{R}$ and*

$$\int |y|^j v_k(x) \bar{p}_m^\varepsilon(x, y)\, dy \leq \tilde{q}(x) e^{l_k M_\delta^2 x^2}$$

*for all $x \in I_\delta^\varepsilon$ where $q, \tilde{q}$ are polynomials of degree $j$ and $l_k = k/(1 - 2\sigma_{ub}^2 k)$.*

The proof uses simple complete the square arguments along with standard properties of Gaussian kernels and is omitted. Applying Lemma 7 with $j = 0$ yields:

**PROPOSITION 2.** *Let $k < \frac{1}{2\sigma_{ub}^2}$. Then $T_s$ is a bounded linear operator from $W_k$ into $W_{l_k c_s^2}$ and for any $\delta > 0$, $T_m^\varepsilon$ is a bounded linear operator from $W_{I_\delta^\varepsilon, k}$ to $W_{I_\delta^\varepsilon, m(k,\delta)}$ with $l_k$ as in Lemma 7 and $m(k,\delta) := l_k M_\delta^2$.*

**REMARK 4.** Note that Lemmas 6 and 7 (and hence, Proposition 2) do not rely on $|c_s| < 1$ and hold equally well if in fact $|c_s| \geq 1$.

With the help of Proposition 2, we can establish a number of useful bounds on growth rates in $W_k$. The key point is that since $|c_s| < 1$, if $k$ is positive, but sufficiently small, we can make $l_k c_s^2 < k$. We use the notation

$$S_{i,j}\phi(x) = \int \phi(y)|x|^i|y|^j p_s(x,y)\,dy$$

and

$$S_{i,j}^\varepsilon\phi(x) = \int \phi(y)|x|^i|y|^j \bar{p}_m^\varepsilon(x,y)\,dy$$

for $\phi \in W_k$ and $x \in \mathbb{R}$.

**COROLLARY 1.** *For any $0 < k < \frac{1-c_s^2}{2\sigma_{ub}^2} =: k_c$ and $n,m \in \mathbb{N}$, there exist positive constants $K_1, L_1 > 0$ depending only on $k, c_s, n, m$ so that*

$$\frac{|S_{i,j}\phi(x)|}{v_k(x)} \leq K_1\|\phi\|_k e^{-L_1 x^2}$$

*$\forall i \leq m, j \leq n,\ x \in \mathbb{R},$ and $\phi \in X$.*

**PROOF.** Let $\phi \in W_k$ and $k' > k$. Then $g\phi \in W_{k'}$ for any polynomial $g$ and $\|g\phi\|_{k'} \leq \|g\|_{k'-k}\|\phi\|_k$ so Proposition 2 implies that

$$\frac{|S_{i,j}\phi(x)|}{v_k(x)} \leq K_1'|x|^i\|\phi\|_k e^{(l_{k'}c_s^2 - k)x^2}$$

for some $K_1' > 0$. Letting $k' \searrow k$ and noting that $0 < k < k_c$ implies $l_k c_s^2 < k$, we obtain the result. $\square$

We also have the analogous result for $S_{i,j}^\varepsilon$. The proof is again a direct consequence of Proposition 2.

COROLLARY 2. *Suppose that $M_\delta < 1$. Then if $0 < k < k_\delta := \frac{1 - M_\delta^2}{1 - 2\sigma_{ub}^2 k}$ and $m, n \in \mathbb{N}$, we can find positive constants $K_2, L_2$ depending only on $k, n, m$ and $M_\delta$ such that*

$$\frac{|S_{i,j}^\varepsilon \phi(x)|}{v_k(x)} \leq K_2 \|\phi\| e^{-L_2 x^2}$$

*$\forall i \leq m, j \leq n, \varepsilon > 0, x \in I_\delta^\varepsilon$, and $\phi \in W_k$.*

Proposition 2 and its corollaries are most useful when $|x|$ is large since the bounds we obtain are then exponentially small. The purposes of the next proposition will be to control the size of $T_s \phi(x)$ when $|x|$ is small, but the support of $\phi$ does not contain $x$, and to control the size of $T_s^\varepsilon - T_m^\varepsilon$. We prove the result in a general form that will also be useful in Section 6.1. We set $U^\varepsilon := U/\varepsilon$ for $U \subset \mathbb{R}$.

PROPOSITION 3. *Suppose $g \in B(U)$ for some $U \subset \mathbb{R}$ and there exists $V \subset \mathbb{R}$ and $\gamma > 0$ such that $d(g(x), V) > \gamma$ for all $x \in U$. Let $\varepsilon > 0$, $g^\varepsilon(x) = \varepsilon^{-1} g(\varepsilon x)$, and let $Y$ be a $N(g^\varepsilon(x), \sigma_\varepsilon^2(x))$ random variable with $\sigma_\varepsilon \in B(\mathbb{R}, \mathbb{R}^+)$ satisfying $0 < \sigma_{lb} < \sigma_\varepsilon(x) < \sigma_{ub}$ for all $x \in \mathbb{R}$, $\varepsilon > 0$ and some constants $\sigma_{lb}, \sigma_{ub}$. Then there exists a $k_\gamma > 0$ so that if $k \in (0, k_\gamma)$ and $\phi \in W_k$ we have*

$$\mathbb{E}[(\phi \mathbf{1}_{V^\varepsilon})(Y)] \leq \frac{\bar{K} \varepsilon}{\gamma} e^{l_k (g^\varepsilon(x))^2} \|\phi\|_k e^{-\bar{L} \gamma^2 / \varepsilon^2}$$

*for all $x \in U^\varepsilon$ and some constants $\bar{K}, \bar{L} > 0$ depending on $k, \sigma_{ub}$, and $\sigma_{lb}$ with $l_k = k/(1 - 2\sigma_{ub}^2 k)$ as in Proposition 2.*

PROOF. If $k < 1/(2\sigma_{ub}^2)$, $x \in U^\varepsilon$ and $\phi \in W_k$,

$$\mathbb{E}[(\phi \mathbf{1}_{V^\varepsilon})(Y)] \leq K_1 \|\phi\|_k \int_{V^\varepsilon} e^{ky^2} e^{-(y - g^\varepsilon(x))^2 / (2\sigma_{ub}^2)} \, dy$$

$$= K_1 \|\phi\|_k e^{l_k (g^\varepsilon(x))^2} \int_{V^\varepsilon} e^{-(y - g^\varepsilon(x)/(1 - 2\sigma_{ub}^2 k))^2 / (2\bar{\sigma}^2)} \, dy$$

$$= K_2 \|\phi\|_k e^{l_k (g^\varepsilon(x))^2} \mathbb{P}\left( \frac{g^\varepsilon(x)}{1 - 2\sigma_{ub}^2 k} + \bar{\sigma} \chi \in V^\varepsilon \right),$$

where $\bar{\sigma}^2 = \sigma_{ub}^2 / (1 - 2\sigma_{ub}^2 k)$ and $\chi$ is standard normal. But from our assumptions on $g$, we have $d(g^\varepsilon(x), V^\varepsilon) > \gamma/\varepsilon$ so that since $1/(1 - 2\sigma_{ub}^2 k) \to 1$ as $k \searrow 0$, we can choose $k$ small enough so that $d(g^\varepsilon(x)/(1 - 2\sigma_{ub}^2 k), V^\varepsilon) > \gamma/(2\varepsilon)$. Therefore,

$$\mathbb{E}[(\phi \mathbf{1}_{V^\varepsilon})(Y)] \leq K_2 \|\phi\|_k e^{l_k (g^\varepsilon(x))^2} \mathbb{P}(\bar{\sigma}|\chi| > \gamma/(2\varepsilon))$$

and the result follows from standard normal distribution tail estimates. $\quad\square$

Using (ii) in Proposition 1 to apply the last result with $g(x) = c_s x$, $\sigma_\varepsilon = \sigma$, $U = V_3$, $V = V_3^c$ and $\gamma = \eta$ yields:

COROLLARY 3. *For all $k$ sufficiently small, $\exists K_3, L_3 > 0$ depending on $c_s, \delta_s, \eta$, and $k$ such that*

$$\frac{|T_s(\phi \mathbf{1}_{(V_3^\varepsilon)^c})(x)|}{v_k(x)} \leq \varepsilon K_3 \|\phi\| e^{-L_3/\varepsilon^2}$$

$\forall x \in V_3^\varepsilon$, *and $\phi \in W_k$.*

At last, we are ready to conclude this section with the desired proof.

PROOF OF THEOREM 5. Without loss of generality, we assume that $M_{\delta_s} < 1$. Choose $k_s$ small enough to apply Corollaries 1–3, let $k < k_s$ and write

$$T_m \phi(x) = \mathbf{1}_{V_3^\varepsilon}(x) \int_{V_3^\varepsilon} \phi(y) p_s(x, y) \, dy.$$

(It is true that the right-hand side of this equation depends on $\varepsilon$ through the cutoff functions, but we do not include this in our notation to make it clear that $T_m$ is related to the limit $T_s$.) We first show that $T_m^\varepsilon = T_s + O(\varepsilon)$. To this end, write

$$T_m^\varepsilon - T_s = (T_m^\varepsilon - T_m) + (T_m - T_s).$$

Then Lemma 6 implies that $p_m^\varepsilon(x, y) = p_s(x, y) + \varepsilon R^\varepsilon(x, y)$ with the remainder term $R^\varepsilon$ satisfying the inequality

$$|R^\varepsilon(x, y)| \leq g_r^\varepsilon(x, y)(\bar{p}_m^\varepsilon(x, y) + p_s(x, y))$$

for any $x, y \in \mathbb{R}$ where $g_r^\varepsilon$ a polynomial in $x, y$ and $\varepsilon$. Therefore, using Corollaries 1 and 2 to control $\int_{V_3^\varepsilon} \phi(y) R^\varepsilon(x, y) \, dy$, we have

$$(14) \qquad \frac{|(T_m^\varepsilon - T_m)(\phi)(x)|}{v_k(x)} \leq O(\varepsilon) \|\phi\|_k$$

for all $\phi \in W_k$, $|x| < \delta_s/\varepsilon$. Furthermore, Corollary 3 implies that

$$(15) \qquad \frac{|T_s(\phi \mathbf{1}_{(V_3^\varepsilon)^c})(x)|}{v_k(x)} \leq O(e^{-L_3/\varepsilon^2}) \|\phi\|_k$$

for all $\phi \in W_k$, $|x| < \delta_s/\varepsilon$ and Corollary 1 yields

$$(16) \qquad \frac{|T_s \phi(x)|}{v_k(x)} \leq O(e^{-L_1 \delta^2/\varepsilon^2}) \|\phi\|_k$$

for all $\phi \in W_k$, $|x| \geq \delta_s/\varepsilon$. Combining (15) and (16), we have

$$\tag{17} \frac{|(T_m - T_s)(\phi)(x)|}{v_k(x)} \leq O(e^{-L/\varepsilon^2}) \|\phi\|_k$$

for all $x \in \mathbb{R}$, $\phi \in W_k$ and some constant $L > 0$. (14) and (17) together imply that $T_m^\varepsilon = T_s + O(\varepsilon)$ in $W_k$ so the proof Theorem 5 will be complete if we can show that

$$\tag{18} \frac{|(T_s^\varepsilon - T_m^\varepsilon)(\phi)(x)|}{v_k(x)} \leq O(e^{-L/\varepsilon^2}) \|\phi\|_k$$

for all $\phi \in W_k$, $x \in \mathbb{R}$. (18) can be obtained by applying Proposition 3 multiple times with $g(x) = f(x)$, $\sigma_\varepsilon^2(x) = \sigma^2(\varepsilon x)$, $U = V_3$, and $V = V_3 + 2\pi n$ for different $n \in \mathbb{Z} \setminus \{0\}$, and noting that $d(f(x), V_3 + 2\pi n) > \eta + 2\pi(|n| - 1)$ whenever $x \in V_3$ (provided $\delta_3 \leq \pi$).  □

Remark 5. Note that we did not use the full force of the expansion for $p_m^\varepsilon$ in Lemma 6; however, we have stated the stronger result anyways to suggest how one can calculate higher order terms in the expansion of $T_s^\varepsilon$. To further explore this idea, we suggest the reader look at the difference between $T_m^\varepsilon$ and $T_s + \varepsilon T_1$ in $W_k$ where

$$T_1\phi(x) := \int \phi(y) g_1(x, y) p_s(x, y) \, dy$$

and $g_1$ is the polynomial appearing in Lemma 6.

4.2. *Expansions of eigenvalues.* With Theorem 5 in hand, we can now apply classical results from Kato's perturbation theory for linear operators (see Chapters 2 and 3 of [6]) to obtain asymptotic expansions for the eigenvalues and eigenvectors of $T_s^\varepsilon$. In order to apply these results, we need to establish the compactness of our operators in $L(W_k)$ and identify the spectrum of $T_s$ is $W_k$. This is done in the following two results.

Lemma 8. $T_s$ and $T_s^\varepsilon, \varepsilon > 0$ are compact operators on $W_k$ for all $k$ sufficiently small.

Proof. An elementary (but somewhat lengthy) calculation using Corollary 2 shows that if $k < k_{\delta_s}$, then $\{T_s^\varepsilon \phi_n(x)\}$ is equicontinuous and equibounded for any sequence of functions $\phi_n \in W_k$ with $\|\phi_n\|_k = 1$ and all $x \in (-\delta_s/\varepsilon, \delta_s/\varepsilon)$. Therefore, by the Arzela–Ascoli theorem, there exists a subsequence $n_j$ and a continuous function $\phi$ defined on $(-\delta_s/\varepsilon, \delta_s/\varepsilon)$ such that $T_s^\varepsilon \phi_{n_j} \to \phi$ uniformly on $[-\delta_s/\varepsilon, \delta_s/\varepsilon]$. If we extend $\phi$ to a function defined on all of $\mathbb{R}$ by setting $\phi(x) = 0$ for $|x| > \delta/\varepsilon$, then $T_s^\varepsilon \phi_{n_j} \to \phi$ in $W_k$ as $n_j \to \infty$ since $T_s^\varepsilon \phi(x) = 0$ if $|x| > \delta_s/\varepsilon$. Therefore, $T_s^\varepsilon$ is compact, for $\varepsilon > 0$. A similar argument applies to $T_s$ with the additional use of Corollary 1 in the final step to control the size of $T_s \phi_{n_j}(x)$ for large $x$.  □

LEMMA 9. $\sigma_{W_k}(T_s) = \{c_s^n\}_{n \geq 0} \cup \{0\}$ *for any $k$ sufficiently small. Furthermore, each $c_s^n$ is a simple eigenvalue with corresponding eigenfunctions $\phi_{s,n}$ and eigenmeasures $\phi_{s,n}^*$.*

PROOF. Follows from Lemmas 5 and 8, and the fact that $W_k \subset L^2(\mu)$ for all $k < \alpha^2$. □

In the end our hard work pays off and we can finish with the following:

PROOF OF THEOREM 4. Theorem 5 and Remark 3 along with the classical results on perturbation theory for linear operators (see [6], Chapters 2 and 3) and Lemma 9 imply that for any given $r > 0$ and $k$ sufficiently small (without loss of generality, we can assume $k < \alpha^2$), $\exists \varepsilon_{s,r}, L_{s,r}, K_{s,r} > 0$ so that for all $\varepsilon < \varepsilon_{s,r}$, any eigenvalue of $T_s^\varepsilon$ in $W_k$ with modulus greater than $r$ is a simple eigenvalue of the form $\lambda_{s,j}^\varepsilon = c_s^j + \varepsilon \lambda_{s,j,1}^\varepsilon$ for some $j \geq 0$ with $|\lambda_{s,j,1}^\varepsilon| \leq L_{s,r}$ and the corresponding eigendensities are multiples of $\phi_{s,j}^\varepsilon(x) := (h_j(\alpha x) + \varepsilon \psi_{s,j}^\varepsilon(x)) \mathbf{1}_{V_3^\varepsilon}(x)$ with $\|\psi_{s,j}^\varepsilon\|_{k-\alpha^2} \leq K_{s,r}$. Assuming without loss of generality that $k < \alpha^2$ yields the appropriate bound for $\psi_{s,j}^\varepsilon$. We have already argued that the nonzero eigenvalues of $T_s^\varepsilon$ in $W_k$ are the same as the nonzero eigenvalues of $T_{33}^\varepsilon$ in $B(V_3)$ and the eigendensities for $T_{33}^\varepsilon$ can be obtained from the eigendensities of $T_s$ by applying $U^\varepsilon$. After re-identifying $[-\pi/2, \pi/2]$ with $S^1$ and 0 with $x_s$, we obtain the result. □

This completes our analysis of the operator near the stable fixed point and explains the first set of eigenvalues appearing in (i) of Theorem 3. We now move to the study of the operator in a neighborhood of the unstable fixed point, which thankfully turns out to be essentially the same.

## 5. The local story near an unstable fixed point.
Our main result for this case is:

THEOREM 6. *For any $r > 0$, $\exists \varepsilon_{u,r}, L_{u,r} > 0$ so that $\forall \varepsilon < \varepsilon_{u,r}$, any eigenvalue of $T_{11}^\varepsilon$ in $B(V_1)$ with modulus greater than $r$ is a simple eigenvalue of the form*

$$\lambda_{u,j}^\varepsilon = |c_u|^{-1} c_u^{-j} + \varepsilon \lambda_{u,j,1}^\varepsilon$$

*for some $j \geq 0$ with $|\lambda_{u,j,1}^\varepsilon| \leq L_{u,r}$, $\forall j \geq 0, \varepsilon < \varepsilon_{u,r}$. Furthermore, $\exists K_{u,r} > 0$ such that the eigenfunctions of $T_{11}^\varepsilon$ corresponding to $\lambda_{u,j}^\varepsilon$ are multiples of*

$$\left[ h_j\left( \frac{\beta(x - x_u)}{\varepsilon} \right) + \varepsilon \psi_{u,j}^\varepsilon\left( \frac{x - x_u}{\varepsilon} \right) \right] \mathbf{1}_{V_1}(x)$$

*with $h_j, \beta$ as in Theorem 3 and*

$$\sup_{x \in \mathbb{R}} (|\psi_{u,j}^{\varepsilon}(x)|e^{kx^2}) \leq K_{u,r}$$

*for all $\varepsilon < \varepsilon_{u,r}$ and some $k > 0$.*

To motivate the proof, we identify $S^1$ with $[-\pi/2, \pi/2]$ by identifying $x_u$ with $0$ (so that $V_1 = (-\delta_u, \delta_u)$, $f(0) = 0$ and $c_u = f'(0) \notin [-1, 1]$) and as in Section 4 consider the chain $Y_n^{\varepsilon} = X_n^{\varepsilon}/\varepsilon$. The limit as $\varepsilon \to 0$ is now the transient chain

$$Y_{n+1} = c_u Y_n + \sigma_0 \chi_n$$

and if we let $T_u$ denote the corresponding transition operator, a simple complete the square calculation reveals that

$$T_u(v_{\beta^2}\phi)(x) = \frac{1}{|c_u|} v_{\beta^2}(x)\mathbb{E}[\phi(x/c_u + (\sigma_0/c_u)\chi)]$$

$\forall \phi \in B$ if $\beta = \sqrt{(c_u^2 - 1)/2\sigma_0}$. Therefore,

$$\frac{|c_u|}{v_{\beta^2}(x)} T_u(v_{\beta^2}\phi)(x) = \mathbb{E}[\phi(x/c_u + (\sigma_0/c_u)\chi)]$$

$\forall \phi \in B, x \in \mathbb{R}$. Notice that the right-hand side is the transition operator for the autoregressive chain

$$Z_{n+1} = c_u^{-1} Z_n + (\sigma_0/c_u)\chi_n$$

and Lemma 9 implies its eigenvalues are $c_u^{-n}$ with corresponding eigenfunctions $\phi_n(x) = H_n(\sqrt{(1 - c_u^{-1})/(2\sigma_0 c_u^{-1})}x) = H_n(\beta x)$. Therefore, $T_u$ will have eigenvalues $|c_u|^{-1}c_u^{-n}$ with corresponding eigenfunctions $\phi_n(x)/v_{\beta^2}(x) = h_n(\beta x)$ as required.

The argument that these are the appropriate limits for the eigenvalues/eigenfunctions of $T_{11}^{\varepsilon}$ can be made rigorous by following the arguments used in Section 4. To this end, we define $T_{11}^{\varepsilon} : W_{k,V_1} \to W_{k,V_1}$, rescale as in Section 4, and extend the resulting operator to $W_k$. We call the extended operator $T_u^{\varepsilon}$ to make clear the analogy with the stable case. From Proposition 2 (see also Remark 4) $T_u^{\varepsilon}$ maps $W_{V_1^{\varepsilon},k}$ to $W_{V_1^{\varepsilon},m(k,\delta_u)}$ for all $k$ but since $|c_u| > 1$, we can check that $m(k, \delta_u) < k$ only if $k < 0$. From our work in Section 4, this suggests that we should look at the limiting behavior of $T_u^{\varepsilon}$ in $W_k$ for some $k < 0$. The next theorem shows that this suggestion is a good one.

THEOREM 7.   *There exists $k_u < 0$ such that $T_u^{\varepsilon} = T_u + O(\varepsilon)$ in $L(W_k)$ for all $k \in (k_u, 0)$.*

PROOF. Consider the operators $\hat{T}_u^\varepsilon := |c_u| V^{-1} \circ T_u^\varepsilon \circ V$ and $\hat{T}_u := |c_u| V^{-1} \circ T_u \circ V$ where $V : W_{k-\beta^2} \to W_k$ is defined by $V\phi(x) = v_{\beta^2}(x)\phi(x)$. From the preceding dialogue, we have

$$\hat{T}_u\phi(x) = \int \phi(y)\hat{p}_u(x,y)$$

for all $\phi \in W_k$ where

$$\hat{p}_u(x,y) = \frac{1}{\sqrt{2\pi}\sigma_0 |c_u|^{-1}} e^{-(y-c_u^{-1}x)^2/(2\sigma_0^2 c_u^{-2})}$$

and a similar calculation reveals that

$$\hat{T}_u^\varepsilon\phi(x) = \mathbf{1}_{V_1^\varepsilon}(x) h^\varepsilon(x) \mathbb{E}[(\phi \mathbf{1}_{V_1^\varepsilon})(F_u^\varepsilon(x) + c_u^{-1}\sigma(\varepsilon x)\chi)],$$

$$= \mathbf{1}_{V_1^\varepsilon}(x) \int_{V_1^\varepsilon} \phi(y) h^\varepsilon(x) \hat{p}_u^\varepsilon(x,y)\, dy,$$

where $F_u^\varepsilon(x) = f_u^\varepsilon(x)/c_u^2$, $h^\varepsilon(x) = e^{\beta^2[x^2 - (f_u^\varepsilon(x)/c_u)^2]}$, and

$$\hat{p}_u^\varepsilon(x,y) = \sum_{n\in\mathbb{Z}} \frac{c_u}{\sqrt{2\pi}\sigma(\varepsilon x)} e^{-c_u^2(y + 2\pi n/\varepsilon - F_u^\varepsilon(x))^2/(2\sigma^2(\varepsilon x))}.$$

Since $(F_u^\varepsilon)'(0) = c_u^{-2}(f_u^\varepsilon)'(0) = c_u^{-1} \in (-1,1)$, and $h^\varepsilon(x) \approx 1$ for small $x$, $\tilde{p}_u^\varepsilon$ has the same basic form as the transition density $\varepsilon\tilde{p}(\varepsilon\cdot, \varepsilon\cdot)$ for the operator $T_s^\varepsilon$ discussed in the previous section. We leave it to the reader to check if $k'$ is sufficiently small and positive, (14)–(18) from the proof of Theorem 5 hold if we replace the pair $T_s^\varepsilon, T_s$ with $\hat{T}_u^\varepsilon, \hat{T}_u$ (and $p_m^\varepsilon$ with the corresponding main part of $\hat{p}_u^\varepsilon$) which implies that $\hat{T}_u^\varepsilon = \hat{T}_u + O(\varepsilon)$ in $W_{k'}$ for $k' > 0$ small. Applying $V$ and $V^{-1}$ to this equation yields the result with $k_u = k' - \beta^2$. □

Theorem 6 then follows from Kato's perturbation theory (see the proof of Theorem 4).

## 6. General periodic orbits.
Having completed our analysis in the case $f$ has only two fixed points, we move on to discuss the issues involved in dealing with general periodic orbits.

6.1. *Stable period two orbit.* We consider the behavior of (3) when $f$ again has two fixed points $x_s, x_u$ but now $f'(x_s), f'(x_u) \notin [-1,1]$ so that $x_s$ and $x_u$ are both unstable. In addition to these two fixed points, we suppose $f$ also has a stable period two orbit $P = \{x_1, x_2\}$ with $(f^2)'(x_1) = (f^2)'(x_2) = f'(x_1)f'(x_2) \in (-1,1)$ and assume that all orbits of $x_{n+1} = f(x_n)$ converge to $P$ if $x_0 \notin \{x_s, x_u\}$. The reader will notice in Theorem 8 below that the eigenvalues in (ii) below also appeared in Theorem 3 as contributions from the unstable fixed point(s). The eigenvalues in (i) are the new contributions from the period two orbit.

THEOREM 8.   *For any $r > 0$, we can decompose $T^\varepsilon = T_{lp}^\varepsilon + T_{up}^\varepsilon$ so that $\forall \varepsilon$ sufficiently small, $\|T_{lp}^\varepsilon\|_\infty < r$ and any eigenvalue of $T_{up}^\varepsilon$ with modulus greater than $r$ is of the form $\lambda^\varepsilon = \lambda + O(\varepsilon)$ where:*

(i)  $\lambda = \sqrt{(c_1 c_2)^j}$ *for some $j \geq 0$ and some branch of $\sqrt{\cdot}$*

*or*

(ii)  $\lambda = |c_k|^{-1} c_k^{-j}$ *for some $j \geq 0$ with $k = s$ or $k = u$.*

*The eigenfunctions corresponding to the limiting eigenvalue in* (ii) *are of the form*

$$a_1 h_n\left(\frac{\beta(x - x_k)}{\varepsilon}\right) + O(\varepsilon)$$

*for some constant $a_1$ with $\beta = \sqrt{(c_k - 1)^2/(2\sigma^2(x_k))}$, $k = s$ or $k = u$. The eigendensities corresponding to the limiting eigenvalues in* (i) *are of the form*

$$\sum_{i=1}^{2} a_i h_n\left(\frac{\alpha_i(x - x_i)}{\varepsilon}\right) + O(\varepsilon)$$

*for some constants $a_i$ with $\alpha_i = \alpha/\sigma_i$, $i = 1, 2$, $\alpha = \sqrt{(1 - c_1 c_2)^2/2}$, $\sigma_1 = \sqrt{c_2^2 \sigma^2(x_1) + \sigma^2(x_2)}$, and $\sigma_2 = \sqrt{c_1^2 \sigma^2(x_2) + \sigma^2(x_1)}$.*

For more information on the $O(\varepsilon)$ terms in the eigenvector expansions see Theorem 3. Note that the difference in the scaling factors $\alpha_j$ imply that the limiting eigendensities corresponding to (i) will have different spread near $x_1$ and $x_2$ (see Figure 4 in Section 2).

The basic form of the proof of Theorem 8 closely resembles the proof of Theorem 3 and the remainder of this section is dedicated to an outline of the steps involved. The main difference comes in the analysis of the part of $T^\varepsilon$ near the stable period two orbit. First, we have the following analog to Proposition 1 which allows us to break up the circle into regions determined by the different actions of $f$.

PROPOSITION 4.   *There exist neighborhoods $V_1 := B_{\delta_u}(x_u)$, $V_2 := B_{\delta_s}(x_s)$, $U_1 := B_{\delta_1}(x_1)$, $U_2 := B_{\delta_2}(x_2)$, and constants $\eta > 0$, $N \in \mathbb{N}$ such that:*

(i)  $d(f(x), V_i) > \eta$ *for every $x \notin V_i$, $i = 1, 2$.*

(ii)  $f(U_1) \subset U_2$, $f(U_2) \subset U_1$ *with $d(f^2(x), U_i^c) > \eta$, and $d(f(x), U_j^c) > \eta$ for every $x \in U_i$, $i, j = 1, 2$, $j \neq i$.*

(iii)  $f^n(x) \in V_4 =: U_1 \cup U_2$ *for every $x \in V_3 := S^1 \setminus (V_1 \cup V_2 \cup V_4)$ and $n \geq N$.*

The proof is similar to the proof of Proposition 1 and is omitted. This splitting leads to the decomposition $T^\varepsilon = (T_{ij}^\varepsilon)_{i,j=1}^{4}$, with $T_{ij}^\varepsilon : B(V_j) \to B(V_i)$ given by (7) for $i, j = 1, 2, 3, 4$. Proposition 4 then tells us that $T^0$ has a block decomposition of the form

$$T^0 = \begin{bmatrix} T_{11}^0 & 0 & T_{13}^0 & T_{14}^0 \\ 0 & T_{22}^0 & T_{23}^0 & T_{24}^0 \\ 0 & 0 & T_{33}^0 & T_{34}^0 \\ 0 & 0 & 0 & T_{44}^0 \end{bmatrix}$$

with $T_{33}^n = 0$ for all $n \geq N$. Furthermore, from property (ii) of Proposition 4, we can further decompose $T_{44}^0$ with respect to $V_4 = U_1 \cup U_2$ into operators $T_{4ij} : B(U_j) \to B(U_i)$ so that

$$T_{44}^0 = \begin{bmatrix} 0 & T_{412}^0 \\ T_{421}^0 & 0 \end{bmatrix}.$$

Lemma 2 can then be used to show that if $\varepsilon > 0$, all the 0 terms in the above decompositions will be replaced in the corresponding decompositions of $T^\varepsilon$ and $T_{44}^\varepsilon$ by terms that are $O(\varepsilon e^{-K/\varepsilon^2})$ as $\varepsilon \to 0$, yielding the appropriate $T_{lp}^\varepsilon$ and $T_{up}^\varepsilon$ terms in Theorem 8. The existence of the eigenvalues in (ii) and the form of the corresponding eigenfunctions then follows directly from Theorem 6 in Section 5. Therefore, our analysis will be complete once we show the eigenvalues of

$$T_4^\varepsilon := \begin{bmatrix} 0 & T_{412}^\varepsilon \\ T_{421}^\varepsilon & 0 \end{bmatrix}$$

correspond to (i) with the appropriate eigendensities. Note that since this operator has a 0 diagonal, its eigenvalues will be given by $\sqrt{\lambda}$ where $\lambda$ is an eigenvalue of $S^\varepsilon := T_{412}^\varepsilon T_{421}^\varepsilon : B(U_1) \to B(U_1)$ and $\sqrt{\cdot}$ denotes the multi-valued complex root function and the corresponding eigendensities will be linear combinations of the eigendensities for $S^\varepsilon$ and $\tilde{S}^\varepsilon = T_{421}^\varepsilon T_{412}^\varepsilon : B(U_2) \to B(U_2)$.

To determine the spectrum of $S^\varepsilon$, we map $S^1$ to $[-\pi/2, \pi/2]$ in such a way that $x_1 \leftrightarrow 0$ [and hence, $x_2 \leftrightarrow f(0)$ with $f(x_2) = 0$]. By definition, $S^\varepsilon \phi(x)$ is zero unless $X_0^\varepsilon \in U_1$, $X_1^\varepsilon \in U_2$, and $X_2^\varepsilon \in U_1$ so to calculate its spectrum, we alternate our re-scaling and look at the chain $Y_n^\varepsilon$ defined by: $Y_0^\varepsilon := X_0^\varepsilon / \varepsilon$,

$$Y_{2n-1}^\varepsilon := \frac{X_{2n-1}^\varepsilon - f(0)}{\varepsilon} = f_1^\varepsilon(Y_{2(n-1)}^\varepsilon) + \sigma_1^\varepsilon(Y_{2(n-1)}^\varepsilon)\chi_0,$$

$$Y_{2n}^\varepsilon := \frac{X_{2n}^\varepsilon}{\varepsilon} = f_2^\varepsilon(Y_{2n-1}^\varepsilon) + \sigma_2^\varepsilon(Y_{2n-1}^\varepsilon)\chi_1$$

with $f_1^\varepsilon(x) := \varepsilon^{-1}(f(\varepsilon x) - f(0)) \to f'(0)x$, $f_2^\varepsilon(x) := \varepsilon^{-1}f(\varepsilon x + f(0)) \to f'(f(0))x$, $\sigma_1^\varepsilon := \sigma(\varepsilon x) \to \sigma(0)$, and $\sigma_2^\varepsilon(x) := \sigma(\varepsilon x + f(0)) \to \sigma(f(0))$ as $\varepsilon \to 0$. Therefore, if $\varepsilon$ is small, the corresponding re-scaled version of $S^\varepsilon$ should be close to the two-step transition operator for the linear chain

$$Y_{2n-1} := c_1 Y_{2(n-1)} + \sigma(0)\chi_0,$$

$$Y_{2n} := c_2 Y_{2n-1} + \sigma(f(0))\chi_1,$$

where $c_1 = f'(0)$ and $c_2 = f'(f(0))$. But

$$Y_{2n+2} = c_1 c_2 Y_{2n} + (c_2 \sigma(0)\chi_{2n} + \sigma(f(0))\chi_{2n+1}) \stackrel{d}{=} c_1 c_2 Y_{2n} + \sigma_1 \tilde{\chi}_{2n}$$

with $\sigma_1^2 = c_2^2 \sigma^2(0) + \sigma^2(f(0))$ and $\tilde{\chi}_{2n}$ a family of i.i.d. standard normal random variables, which we recognize as the autoregressive scheme previously encountered in Section 4. Since $|c_1 c_2| < 1$, Lemma 9 implies the transition operator for this chain has eigenvalues $(c_1 c_2)^n$, $n \geq 0$. Therefore, the eigenvalues of $S^\varepsilon$ should also be close to $(c_1 c_2)^n$ for small $\varepsilon$ yielding the eigenvalues (ii) in Theorem 8.

To make this argument rigorous, we again make use of the weighted sup-norm spaces $W_k$ defined in Section 4. In the next result, we shall use $S^\varepsilon$ to denote the two-step transition operator for the re-scaled chain $Y_n^\varepsilon$ and $S$ to denote the two step transition operator of the limiting chain $Y_n$.

THEOREM 9. *For any sufficiently small $k > 0$, we have $S^\varepsilon = S + O(\varepsilon)$ in $W_k$.*

PROOF. We assume $\sigma$ is constant for notational purposes and leave the general case to the reader. Write $p_1^\varepsilon(x,y) := \varepsilon p^\varepsilon(\varepsilon x, \varepsilon y + f(0)), p_2^\varepsilon(x,y) := \varepsilon p^\varepsilon(\varepsilon x + f(0), \varepsilon y)$ for the main parts of the transition densities for $Y_{2n-1}^\varepsilon, Y_{2n}^\varepsilon$ and $p_1, p_2$ for the transition densities of $Y_{2n-1}$ and $Y_{2n}$, respectively. Define

$$S_i^\varepsilon \phi(x) := \mathbf{1}_{U_i^\varepsilon}(x) \int_{U_j^\varepsilon} p_i^\varepsilon(x,y)\, dy$$

and

$$S_i \phi(x) := \mathbf{1}_{U_i^\varepsilon}(x) \int_{U_j^\varepsilon} p_i(x,y)\, dy$$

for $i,j = 1,2$, $i \neq j$ where here $U_i^\varepsilon := (-\delta_i/\varepsilon, \delta_i/\varepsilon)$. Finally, let $S_m^\varepsilon := S_2^\varepsilon S_1^\varepsilon$ and $S_m = S_2 S_1$.

As in the proof of Theorem 5, we will show that $S_m^\varepsilon = S + O(\varepsilon)$ by bounding $S_m^\varepsilon - S_m$ and $S_m - S$ and then show that $S^\varepsilon - S_m^\varepsilon$ is small. Our first task will be proving

$$(19) \qquad \frac{|S_m^\varepsilon \phi(x) - S_m \phi(x)|}{v_k(x)} \leq O(\varepsilon)\|\phi\|_k$$

for all $\phi \in W_k, |x| < \delta_1/\varepsilon$, the natural analog of (14).

To prove (19) let $k > 0$, $x \in U_1^\varepsilon$ and $\phi \in W_k$. Then

$$|S_m^\varepsilon \phi(x) - S_m \phi(x)| = \left| \int_{U_1^\varepsilon} \int_{U_2^\varepsilon} \phi(z)(p_1(x,y)p_2(y,z) - p_1^\varepsilon(x,y)p_2^\varepsilon(y,z))\, dy\, dz \right|.$$

Using

$$p_1(x,y)p_2(y,z) - p_1^\varepsilon(x,y)p_2^\varepsilon(y,z) = (p_1(x,y) - p_1^\varepsilon(x,y))p_2(y,z)$$
$$+ p_1^\varepsilon(x,y)(p_2(y,z) - p_2^\varepsilon(y,z)),$$

we have

$$(20) \qquad |S_m^\varepsilon \phi(x) - S_m \phi(x)| \le |(S_1 - S_1^\varepsilon)[S_2\phi](x)| + |S_1^\varepsilon[(S_2 - S_2^\varepsilon)\phi](x)|.$$

Since $f_1^\varepsilon(y) \to c_1 y$ and $f_2^\varepsilon(y) \to c_2 y$ as $\varepsilon \to 0$, Lemma 6 implies that

$$|p_1(x,y) - p_1^\varepsilon(x,y)| \le O(\varepsilon) g_{r,1}^\varepsilon(x,y)(p_1(x,y) + p_1^\varepsilon(x,y))$$

and

$$|p_2(y,z) - p_2^\varepsilon(y,z)| \le O(\varepsilon) g_{r,2}^\varepsilon(y,z)(p_2(y,z) + p_2^\varepsilon(y,z)),$$

where $g_{r,1}^\varepsilon$ and $g_{r,2}^\varepsilon$ are polynomials in $\varepsilon, x, y$ and $\varepsilon, y, z$, respectively. Applying these bounds to (20) and dividing by $v_k(x)$, we have

$$(21) \qquad \begin{aligned} \frac{|S_m^\varepsilon \phi(x) - S_m \phi(x)|}{v_k(x)} &\le O(\varepsilon)\left[\frac{|S_{r,1}(S_2\phi)(x)|}{v_k(x)} + \frac{|S_{r,1}^\varepsilon(S_2\phi)(x)|}{v_k(x)}\right] \\ &\quad + O(\varepsilon)\left[\frac{|S_1^\varepsilon(S_{r,2}\phi)(x)|}{v_k(x)} + \frac{|S_1^\varepsilon(S_{r,2}^\varepsilon\phi)(x)|}{v_k(x)}\right], \end{aligned}$$

where the $S_{r,i}^\varepsilon$ are defined as

$$S_{r,i}^\varepsilon \phi(x) := \int_{U_j^\varepsilon} \phi(y) g_{r,i}(x,y) p_i^\varepsilon(x,y)\, dy$$

for $i, j = 1, 2$, $j \ne i$ and similarly for $S_{r,i}$, $i = 1, 2$.

To bound the four integrals on the right of (21), we appeal to Proposition 2 in Section 4. For instance, $S_2 : W_k \to W_{l_1}$ with $l_1 = l_1(k) = c_2^2 k/(1 - 2\sigma^2 k)$ and $S_1 : W_{l_1} \to W_{l_2}$ where $l_2 = c_1^2 l_1/(1 - 2\sigma^2 l_1)$ so that $S_1 S_2 : W_k \to W_{l_k}$ where

$$l_k = c^2 k/[(1 - 2\sigma^2 k)(1 - 2\sigma^2 l_1)] = c^2 k/(1 - 2\sigma_1^2 k) > 0$$

provided $0 < k < (1 - c^2)/(2\sigma_1^2)$. Furthermore, since $|c| < 1$, we have $l_k < k$ for this same range of $k$ values so that we can apply the argument used in Corollary 1 to bound the first term on the right-hand side of (21). Bounds for the other terms follows in a similar manner, completing the proof of (19).

Our next task is to look at the difference between $S_m$ and $S$. Since $S$ has the same form as $T_s$ in Section 4 (with $c = c_1 c_2$ replacing $c_s$), Proposition 2 and the argument from Corollary 1 imply that if $k < (1 - c^2)/(2\sigma_1^2)$ and $|x| \ge \delta_1/\varepsilon$,

$$(22) \qquad \frac{|(S_m - S)\phi(x)|}{v_k(x)} = \frac{|S\phi(x)|}{v_k(x)} \le O(e^{-L_1 \delta_1^2/\varepsilon^2}) \|\phi\|_k$$

for all $\phi \in W_k$. Furthermore, if $|x| < \delta_1/\varepsilon$, we can write

$$
\begin{aligned}
(23) \quad |(S - S_m)(\phi)(x)| &= |\mathbb{E}^x[\phi(Y_2)(1 - \mathbf{1}_{U_1^\varepsilon}(Y_2)\mathbf{1}_{U_2^\varepsilon}(Y_1))]| \\
&\leq \mathbb{E}^x[|\phi(Y_2)|(1 - \mathbf{1}_{U_2^\varepsilon}(Y_1))] \\
&\quad + \mathbb{E}^x[|\phi(Y_2)|(1 - \mathbf{1}_{U_1^\varepsilon}(Y_2))\mathbf{1}_{U_2^\varepsilon}(Y_1)] \\
&= \mathbb{E}^x[\mathbb{E}^{Y_1}[|\phi(Y_2)|](1 - \mathbf{1}_{U_2^\varepsilon}(Y_1))] \\
&\quad + \mathbb{E}^x[\mathbb{E}^{Y_1}[|\phi(Y_2)|(1 - \mathbf{1}_{U_1^\varepsilon}(Y_2))]\mathbf{1}_{U_2^\varepsilon}(Y_1)].
\end{aligned}
$$

Now from (ii) in Proposition 4, we can assume without loss of generality that $d(c_1 x, I_{\delta_2}^c) > \eta$ if $x \in I_{\delta_1}$ where as before $I_\delta = (-\delta, \delta)$. Therefore, to control the first term on the right of (23), we use Lemma 7 to bound $E^{Y_1}[\phi(Y_2)]$ and then apply Proposition 3 with $g(x) = c_1 x$, $\sigma_\varepsilon = 1$, $U = I_{\delta_1}$, and $V = I_{\delta_2}^c$. Similarly, to control the second term on the right of (23), we apply Proposition 3 with $g(x) = c_2 x$, $\sigma_\varepsilon = 1$, $U = I_{\delta_2}$, and $V = I_{\delta_1}^c$ and then use Lemma 7 to bound the expectation. These bounds for the two terms on the right of (23) yield the inequality

$$
(24) \qquad \frac{|S\phi(x) - S_m\phi(x)|}{v_k(x)} \leq O(e^{-L/\varepsilon^2})\|\phi\|_k
$$

for all $\phi \in W_k, |x| < \delta_1/\varepsilon$ and some $L > 0$ provided $k$ is sufficiently small. (19), (22) and (24) complete that proof that $S_m^\varepsilon = S + O(\varepsilon)$. We leave it to the reader to check that

$$
\frac{|S^\varepsilon\phi(x) - S_m^\varepsilon\phi(x)|}{v_k(x)} \leq O(e^{-L/\varepsilon^2})\|\phi\|_k
$$

for $\phi \in W_k$ and $x \in \mathbb{R}$ as well, which then yields Theorem 9. $\square$

Theorem 9 (along with the perturbation theory arguments used to derive Theorem 4) tells us that for small $\varepsilon$, the top eigenvalues of $S^\varepsilon$ will be $(c_1 c_2)^n + O(\varepsilon)$ and the corresponding eigendensities will be $h_n(\frac{\alpha(x-x_1)}{\varepsilon\sigma_1}) + O(\varepsilon)$. [We have absorbed the cut-off functions that appear in the eigendensity formulas from Theorem 4 in the $O(\varepsilon)$ term since $h_n(\frac{\alpha(x-x_1)}{\varepsilon\sigma_1})$ is concentrated near $x_1$ for small $\varepsilon$ anyways.] We can then apply the off-diagonal structure of $T_4^\varepsilon$ to yield the appropriate limiting eigenvalues and corresponding eigendensities (see the discussion following Proposition 4).

6.2. *Notes on the general case.* The general case described in Section 2 can be handled in the same way as the specific cases we have dealt with in Sections 3 and 6.1 although the details are more tedious. We conclude this paper by remarking on some of the differences.

- The starting point is, as before, to split the circle into regions describing the different actions of $f$. The notation is more complicated, but the end result is the same: we can split the circle into neighborhoods of the different periodic orbits for $f$ and label sets in such a way that $T^\varepsilon$ has an "almost" upper triangular decomposition with respect to the splitting (see Propositions 1 and 4).

- We already know how to deal with the local behavior of the operator near fixed points and stable period two orbits. For stable periodic orbits $P$ of period $p > 2$, we simply note that for small $\varepsilon$, the block corresponding to $P$ in the decomposition of $T^\varepsilon$ is approximately of the form:

$$
\begin{bmatrix}
0 & T_{12}^\varepsilon & 0 & \cdots & 0 \\
\vdots & 0 & T_{23}^\varepsilon & \ddots & \vdots \\
& & \ddots & \ddots & 0 \\
0 & \cdots & & 0 & T_{(p-1)p}^\varepsilon \\
T_{p1}^\varepsilon & 0 & \cdots & & 0
\end{bmatrix},
$$

when $\varepsilon$ is small. One can readily check that any such operator has eigenvalues $\lambda^{1/p}$ where $\lambda$ is an eigenvalue of $T_{12}^\varepsilon T_{23}^\varepsilon \cdots T_{(p-1)p}^\varepsilon T_{p1}^\varepsilon$. A similar argument to the one in Section 6.1 can be used to show that this $p$-step chain has eigenvalues close to $c^n$ where $c$ is equal to the derivative of $f$ along $P$ which yields the eigenvalues $(c^n)^{1/p}$ as desired.

- For an unstable period two orbit $Q$, we identify $S^1$ with $\mathbb{R}$ so that $Q = \{0, f(0)\}$ and again apply the re-scaling argument from Section 6.1 to the block of $T^\varepsilon$ corresponding to $Q$. The limiting two step chain is of the form $Y_{n+2} = cY_n + \sigma_1\chi$ where $c = f'(0)f'(f(0))$ is the derivative along the period two orbit and $\sigma_1^2 = c_2^2\sigma^2(0) + \sigma^2(f(0))$ with $c_2 = f'(f(0))$. Since $|c| > 1$, our work in Section 5 implies that this chain has eigenvalues $|c|^{-1}c^{-n}$ and hence, it seems reasonable to believe that the block of $T^\varepsilon$ corresponding to $Q$ has eigenvalues near $(|c|^{-1}c^{-n})^{1/2}$. To prove this rigorously, we let $S^\varepsilon$ and $S$ be as in Section 6.1 and prove $S^\varepsilon = S + O(\varepsilon)$ in $W_k$ for some $k < 0$ (see Theorem 7). Note that this will require pre and post-multiplying functions by $v_{\beta^2}$ where $\beta^2 = (c^2 - 1)/(2\sigma_1^2)$. Extensions to the case when $Q$ has period $q > 2$ are similar.

**Acknowledgments.** This paper was put together while I was a graduate student at the University of Southern California and I would especially like to thank my advisor, Professor Peter Baxendale, for his countless hours of help in the development of this project. Without his valuable suggestions and guidance, this work would not have been possible. I would also like to thank Professor Paul Newton for pointing out the connections with pseudospectra.

Department of Mathematics
Cornell University
310 Malott Hall
Ithaca, New York 14853-4201
USA
E-mail: jmayberry@math.cornell.edu